% plain-tex file "finite.tex"
% finished on Feb. 8, 2000 (without introduction) and sent to Jeremy
% 2.-5.3.2000, last section 4 added, PS
% 6.3., further corrections, PS
% 12.3., further corrections, proof of Lemma 3.3 completely revised, PS
% 14.3., Jeremy's introduction added, PS
% 17.3., references [DG] and [GGP] added,
%        hard page breaks on p. 14, 18, 19, PS
% 30.4., sentence in introduction about relation to admissibility
%        definition in [ST] updated, PS

% Definitionsfile "defps.tex vom 13.10.99

%  Initialisation der TeX-Konstanten
%
\magnification1200
\pretolerance=100
\tolerance=200
\hbadness=1000
\vbadness=1000
\linepenalty=10
\hyphenpenalty=50
\exhyphenpenalty=50
\binoppenalty=700
\relpenalty=500
\clubpenalty=5000
\widowpenalty=5000
\displaywidowpenalty=50
\brokenpenalty=100
\predisplaypenalty=7000
\postdisplaypenalty=0
\interlinepenalty=10
\doublehyphendemerits=10000
\finalhyphendemerits=10000
\adjdemerits=160000
\uchyph=1
\delimiterfactor=901
\hfuzz=0.1pt
\vfuzz=0.1pt
\overfullrule=5pt
\hsize=146 true mm
\vsize=8.9 true in
\maxdepth=4pt
\delimitershortfall=.5pt
\nulldelimiterspace=1.2pt
\scriptspace=.5pt
\normallineskiplimit=.5pt
\mathsurround=0pt
\parindent=20pt
\catcode`\_=11
\catcode`\_=8
\normalbaselineskip=12pt
\normallineskip=1pt plus .5 pt minus .5 pt
\parskip=6pt plus 3pt minus 3pt
\abovedisplayskip = 12pt plus 5pt minus 5pt
\abovedisplayshortskip = 1pt plus 4pt
\belowdisplayskip = 12pt plus 5pt minus 5pt
\belowdisplayshortskip = 7pt plus 5pt
\normalbaselines
\smallskipamount=\parskip
 \medskipamount=2\parskip
 \bigskipamount=3\parskip
\jot=3pt
%
% Definitionen
%
\def\ref#1{\par\noindent\hangindent2\parindent
 \hbox to 2\parindent{#1\hfil}\ignorespaces}
%
%  fonts
%
        % for running heads
                     % (standard)
   % for unnumbered section titles
   % for   numbered section titles
   % for chapter titles
                   % for authors
\font\tenss=cmss10
\font\sevenss=cmss8 at 7pt
\font\fivess=cmss8 at 5pt
\newfam\ssfam %
\textfont\ssfam=\tenss
\scriptfont\ssfam=\sevenss
\scriptscriptfont\ssfam=\fivess
%
%  Definition der Suffixbuchstaben
%
%
%  Version vom 3.2.1989
%
%
%  erzeugt Buchstaben aus den Standardfonts
%
%
%  (BCDEFGHIKLMNOPQRUZ) : Doppelstrich-rm-Buchstaben
%                         fuer jeden Modus
%
%  (A-Z)(rm|scr|sl|bf|tt|ss) und (a-z)(rm|sl|bf|tt|ss) mathchars
%  ohne hss und vss (!)
%
\catcode`\_=11
\def\suf_fix{}
\def\scaled_rm_box#1{%
 \relax
 \ifmmode
   \mathchoice
    {\hbox{\tenrm #1}}%
    {\hbox{\tenrm #1}}%
    {\hbox{\sevenrm #1}}%
    {\hbox{\fiverm #1}}%
 \else
  \hbox{\tenrm #1}%
 \fi}
\def\suf_fix_def#1#2{\expandafter\def\csname#1\suf_fix\endcsname{#2}}
\def\I_Buchstabe#1#2#3{%
 \suf_fix_def{#1}{\scaled_rm_box{I\hskip-0.#2#3em #1}}
}
\def\rule_Buchstabe#1#2#3#4{%
 \suf_fix_def{#1}{%
  \scaled_rm_box{%
   \hbox{%
    #1%
    \hskip-0.#2em%
    \lower-0.#3ex\hbox{\vrule height1.#4ex width0.07em }%
   }%
   \hskip0.50em%
  }%
 }%
}
\I_Buchstabe B22
\rule_Buchstabe C51{34}
\I_Buchstabe D22
\I_Buchstabe E22
\I_Buchstabe F22
\rule_Buchstabe G{525}{081}4
\I_Buchstabe H22
\I_Buchstabe I20
\I_Buchstabe K22
\I_Buchstabe L20
\I_Buchstabe M{20em }{I\hskip-0.35}
\I_Buchstabe N{20em }{I\hskip-0.35}
\rule_Buchstabe O{525}{095}{45}
\I_Buchstabe P20
\rule_Buchstabe Q{525}{097}{47}
\I_Buchstabe R21 %vorher: 23
\rule_Buchstabe U{45}{02}{54}
\suf_fix_def{Z}{\scaled_rm_box{Z\hskip-0.38em Z}}
\catcode`\"=12
\newcount\math_char_code
\def\suf_fix_math_chars_def#1{%
 \ifcat#1A
  \expandafter\math_char_code\expandafter=\suf_fix_fam
  \multiply\math_char_code by 256
  \advance\math_char_code by `#1
  \expandafter\mathchardef\csname#1\suf_fix\endcsname=\math_char_code
  \let\next=\suf_fix_math_chars_def
 \else
  \let\next=\relax
 \fi
 \next}
%
% \font_fam_suf_fix family suffix buchstaben
%
% definiert \(Buchstabe)(suffix) als mathchar 0(family)`(Buchstabe)
%
%
\def\font_fam_suf_fix#1#2 #3 {%
 \def\suf_fix{#2}
 \def\suf_fix_fam{#1}
 \suf_fix_math_chars_def #3.
}
\font_fam_suf_fix
 0rm
 ABCDEFGHIJKLMNOPQRSTUVWXYZabcdefghijklmnopqrstuvwxyz
\font_fam_suf_fix
 2scr
 ABCDEFGHIJKLMNOPQRSTUVWXYZ
\font_fam_suf_fix
 \slfam sl
 ABCDEFGHIJKLMNOPQRSTUVWXYZabcdefghijklmnopqrstuvwxyz
\font_fam_suf_fix
 \bffam bf
 ABCDEFGHIJKLMNOPQRSTUVWXYZabcdefghijklmnopqrstuvwxyz
\font_fam_suf_fix
 \ttfam tt
 ABCDEFGHIJKLMNOPQRSTUVWXYZabcdefghijklmnopqrstuvwxyz
\font_fam_suf_fix
 \ssfam
 ss
 ABCDEFGHIJKLMNOPQRSTUVWXYZabcdefgijklmnopqrstuwxyz
\catcode`\_=8
\def\Cdss{{\fam\ssfam
    \mkern 4.2 mu \mathchoice%
    {\vrule height 6.5pt depth -.55pt width 1pt}%
    {\vrule height 6.5pt depth -.57pt width 1pt}%
    {\vrule height 4.55pt depth -.28pt width .8pt}%
    {\vrule height 3.25pt depth -.19pt width .6pt}%
    \mkern -6.3mu C}}%
\def\Ndss{{\fam\ssfam I\mkern -2.5mu N}}%
\def\Qdss{{\fam\ssfam
    \mkern 3.8 mu \mathchoice%
    {\vrule height 6.5pt depth -.67pt width 1pt}%
    {\vrule height 6.5pt depth -.7pt width 1pt}%
    {\vrule height 4.55pt depth -.44pt width .7pt}%
    {\vrule height 3.25pt depth -.3pt width .5pt}%
    \mkern -5.9mu Q}}%
\def\Zdss{{\fam\ssfam Z\mkern-8.1mu Z}}%
%
%
%
% Benutzung der Euler-Fraktur-Buchstaben unter TeX
%
\font\teneuf=eufm10 % Euler-Fraktur
%\font\preloaded=eufm9
%\font\preloaded=eufm8
\font\seveneuf=eufm7
%\font\preloaded=eufm6
\font\fiveeuf=eufm5
\newfam\euffam \def\euf{\fam\euffam\teneuf} % \euf is family ?
\textfont\euffam=\teneuf \scriptfont\euffam=\seveneuf
\scriptscriptfont\euffam=\fiveeuf

       \def\gfr{{\euf g}}

       \def\lfr{{\euf l}}

       \def\sfr{{\euf s}}
       
       \def\ufr{{\euf u}}

       \def\xfr{{\euf x}}
       
       \def\zfr{{\euf z}}
\parindent=0pt
\def\dQ{{\Qdss_p}}
\def\dZ{{\Zdss_p}}
\def\Hom{\mathop{\rm Hom}\nolimits}
\def\Ind{{\rm Ind}}
\def\ind{{\rm ind}}
\def\dlongrightarrow{\longrightarrow\hskip-8pt\rightarrow}

\centerline{\bf $U(\gfr)$-finite locally analytic representations}

\medskip

\centerline{\bf P. Schneider, J. Teitelbaum}

\bigskip

In this paper we continue the study of locally analytic
representations of a $p$-adic Lie group $G$ in vector spaces over a
spherically complete non-archimedean field $K$. In [ST], we began with
an algebraic approach to this type of representation theory based on a
duality functor that replaces locally analytic  representations by
certain topological modules over the algebra $D(G,K)$ of locally
analytic distributions. As an application, we established the
topological irreducibility of generic locally analytic principal
series representations of ${\bf GL_{2}}(\dQ)$ by proving the algebraic
simplicity of the corresponding $D({\bf GL_{2}}(\dQ),K)$-modules.

In this paper we further exploit this algebraic point of view. We
introduce a particular category  of ``analytic'' $D(G,K)$-modules
that lie in the image of the duality functor and therefore
correspond to certain locally analytic representations. For
compact groups $G$, these are finitely generated $D(G,K)$-modules
that allow a (necessarily uniquely determined) Fr\'echet topology
for which the $D(G,K)$-action is continuous. For more general
groups, one tests analyticity by considering the action of
$D(H,K)$ for a compact open subgroup $H$ in $G$. The category of
analytic modules has the nice property that any algebraic map
between such modules is automatically continuous. The concept of
analytic module is dual to the concept of strongly admissible
$G$-representation introduced in [ST]. The actual definition can
and will be given in a way that avoids any mention of a topology
on the module.

Next, we study the modules dual to  the traditional smooth
representations of Langlands theory. We show that a smooth
representation gives rise, under duality, to an analytic module
precisely when it is ``strongly admissible''; this is a condition on
the multiplicities with which the irreducible representations of a
compact open subgroup of $G$ appear in the representation. In
particular, if $L$ is a finite extension of $\dQ$ and $G$ is the group
of $L$-points of a connected reductive algebraic group over $L$, then
any smooth representation of finite length is strongly admissible.
This is basically a theorem of Harish-Chandra ([HC]) although we must
use in addition results of Vigneras ([Vig]) to deal with some
complications arising from the fact that we do not assume that our
coefficient field $K$ is algebraically closed.

Given these foundational results, suppose that $G$ is the group of
$L$-points of a split, semisimple, and simply connected  group over
$L$. We completely determine the structure of analytic modules $M$
that are $U(\gfr)$-finite, i.e., that are annihilated by a 2-sided
ideal of finite codimension in the universal enveloping algebra
$U(\gfr)$ of the Lie algebra $\gfr$ of $G$. Such a module can be
decomposed into a finite sum of modules of the form $E\otimes
\Hom(V,K)$ where $E$ is irreducible, finite dimensional, and
algebraic, and $V$ is  smooth and strongly admissible. The dual
representations $E^{\ast}\otimes V$ are irreducible -- in fact, simple
as $K[G]$ modules -- if and only if $V$ is irreducible. Some of the
technical hypotheses on the group $G$ in this section are consequences
of the fact that our coefficient field is not algebraically closed.

We conclude the paper by studying the reducible members of the locally
analytic principal series of ${\bf SL_{2}}(\dQ)$. The corresponding
modules contain a simple submodule such that the quotient is
$U(\gfr)$-finite, and we use our methods to determine the structure of
this quotient. In particular, we obtain the result that the
topological length of the locally analytic principal series is at most
three -- a fact that is due to Morita ([Mor]) by a different method.

\bigskip

{\bf 1. Analytic modules}

\smallskip

We fix fields $\dQ\subseteq L\subseteq K$ such that $L/\dQ$ is
finite and $K$ is spherically complete with respect to a
nonarchimedean absolute value $|\; |$ extending the one on $L$. We
let $G$ be a $d$-dimensional locally $L$-analytic group and
$D(G,K)$ be the corresponding $K$-algebra of $K$-valued
distributions on $G$. Recall ([ST] 2.3) that $D(G,K)$ is an
associative unital $K$-algebra with a natural locally convex
topology in which the multiplication $*$ is separately continuous.
Unless this topology is explicitly mentioned $D(G,K)$ is treated
as an abstract algebra. In the following we want to single out a
certain class of (unital left) $D(G,K)$-modules which seems to
provide a convenient framework for the representation theory of
$G$ over $K$. Let $M$ be a $D(G,K)$-module.

\medskip

{\bf Definition:}

{\it A $K$-linear form $\ell$ on $M$ is called locally analytic
if, for any $m \in M$, the linear form $\lambda \longmapsto \ell
(\lambda m)$ on $D(G,K)$ is continuous.}

\medskip

Clearly the locally analytic linear forms on $M$ form a vector
subspace $M'$ of the full $K$-linear dual $M^{\ast}$ of $M$. We
first consider the case of a compact group $G$. Recall that then
$D(G,K)$ is a $K$-Fr\'echet algebra and as a locally convex
$K$-vector space is reflexive ([ST] 1.1, 2.1, and 2.3).

\medskip

{\bf Definition:}

{\it Suppose $G$ to be compact; a $D(G,K)$-module $M$ is called
analytic if it is finitely generated and if, for any $m \in M$,
there is a locally analytic linear form $\ell$ on $M$ such that
$\ell (m) \neq 0$.}

\medskip

{\bf Proposition 1.1:}

{\it Suppose $G$ to be compact; for a finitely generated
$D(G,K)$-module $M$ the following assertions are equivalent:

i. $M$ is analytic;

ii. $M$ carries a Fr\'echet topology with respect to which it is a
continuous $D(G,K)$-module.}

Proof: We first assume that ii. holds true. Evidently any
continuous linear form on $M$ then is locally analytic. Hence it
follows from the Hahn-Banach theorem that $M$ is analytic. Assume
now vice versa that $M$ is analytic. Choose an epimorphism $\alpha
: D(G,K)^r \dlongrightarrow M$ of $D(G,K)$-modules for some $r
\geq 1$. Then the linear forms $\ell \circ \alpha$ for any
$\ell \in M'$ are continuous and their simultaneous kernel
coincides with the kernel of $\alpha$. In particular the kernel of
$\alpha$ is closed in $D(G,K)^r$ so that the quotient topology via
$\alpha$ on $M$ has the required properties.

\medskip

By the argument in the proof of [ST] 3.5 the above Fr\'echet
topology on an analytic $D(G,K)$-module $M$ is unique and
therefore will be called the {\it canonical\ topology} of $M$. The
continuous dual of $M$ is $M'$ and given the strong topology it is
a vector space of compact type carrying a locally analytic
$G$-representation ([ST] \S\S1 and 3); in particular, the
canonical topology on $M$ is reflexive. Again by [ST] 3.5 any
$D(G,K)$-linear map between two analytic $D(G,K)$-modules is
continuous in the canonical topologies.

\medskip

{\bf Question:} Is any $D(G,K)$-module of finite presentation
analytic ?

\medskip

{\bf Example:} As a consequence of [ST] 4.4 the answer is yes for
the group $G = \dZ$.

\medskip

The above definition of an analytic $D(G,K)$-module for a compact
group is extended to a general group $G$ in the following way.
Note first that for any compact open subgroup $H \subseteq G$ the
algebra $D(H,K)$ is a subalgebra of $D(G,K)$ and that
$$
D(G,K) = \mathop{\bigoplus}\limits_{g \in G/H} \delta_g * D(H,K)\
$$
where $\delta_g$ denotes the Dirac distribution in $g \in G$.

\medskip

{\bf Definition:}

{\it A $D(G,K)$-module $M$ is called analytic if it is analytic as
a $D(H,K)$-module for any compact open subgroup $H \subseteq G$.}

\medskip

{\bf Lemma 1.2:}

{\it Fix a compact open subgroup $H \subseteq G$; a
$D(G,K)$-module $M$ is analytic if it is analytic as a
$D(H,K)$-module.}

Proof: This follows easily from the fact that for any two compact
open subgroups $H$ and $H'$ in $G$ the intersection $H \cap H'$ is
of finite index in $H$ and in $H'$.

\medskip

Suppose that $M$ is an analytic $D(G,K)$-module. One easily checks
that the canonical topology of $M$ as a $D(H,K)$-module is
independent of the choice of the compact open subgroup $H
\subseteq G$, that the $D(G,K)$-action on $M$ is separately
continuous, and that $M'$ is the continuous dual of $M$ and
equipped with the strong topology carries a locally analytic
$G$-representation. Of course, any $D(G,K)$-linear map between two
analytic $D(G,K)$-modules is continuous in the canonical
topologies.

\medskip

{\bf Definition:}

{\it An analytic $D(G,K)$-module is called quasi-simple if it has
no nonzero proper $D(G,K)$-submodules which are closed in the
canonical topology.}

\medskip

An analytic $D(G,K)$-module $M$ is trivially quasi-simple if it is
(algebraically) simple. But, as a consequence of polarity, it also
is quasi-simple (and usually not simple) if $M'$ is a simple
$D(G,K)$-module. For a noncompact $G$ we will see examples of this
later on. We don't know whether such examples also exist for
compact groups.

\bigskip

{\bf 2. Smooth $G$-representations}

\smallskip

In this section we want to see how the smooth representation
theory of $G$ fits into our new framework. We recall that a smooth
$G$-representation $V$ (over $K$) is a $K$-vector space $V$ with a
linear $G$-action such that the stabilizer of each vector in $V$
is open in $G$. (Traditionally one considers smooth
$G$-representations in $\Cdss$-vector spaces; but since the
topology of the coefficient field plays absolutely no role in the
definition this makes a difference only insofar as we do not
require $K$ to be algebraically closed.) Moreover, a smooth
$G$-representation $V$ is called admissible if, for any compact
open subgroup $H \subseteq G$, the vector subspace $V^H$ of
$H$-invariant vectors in $V$ is finite dimensional. Finally,
irreducibility of a smooth representation is always meant in the
algebraic sense.

The unit element in $G$ has a countable fundamental system of open
compact neighborhoods. This implies that the finest locally convex
topology on an admissible $G$-representation $V$ is of compact
type (being the countable locally convex inductive limit of the
finite dimensional spaces $V^H$). Since the orbit maps $\rho_v (g)
:= gv$, for $v \in V$, are locally constant on $G$ we see that any
admissible $G$-representation $V$ becomes a locally analytic
$G$-representation on a vector space of compact type once we equip
$V$ with the finest locally convex topology; as such we denote it
by $V^c$.

Let $\gfr$ denote the Lie algebra of $G$ and let $U(\gfr )$ be the
universal enveloping algebra of $\gfr$. The latter is naturally
included in $D(G,K)$ ([ST] \S2). The action of an $\xfr \in \gfr$
on a locally analytic $G$-representation $W$ is given by
$$
w \to {\xfr}w:={{d}\over{dt}}\exp(t\xfr)w|_{t=0}\leqno{(1)}
$$
where $\exp:\gfr\, {--->}\, G$ denotes the exponential map defined
locally around 0 ([ST] 3.2). In addition Taylor's formula says
that, for each fixed $w\in W$ there is a sufficiently small
neighborhood $U$ of 0 in $\gfr$ such that, for $\xfr\in U$, we
have a convergent expansion
$$
\exp(\xfr)w=\sum_{n=0}^{\infty} {{1}\over{n!}}\xfr^{n}w\
.\leqno{(2)}
$$
The formulas $(1)$ and $(2)$ together imply that the orbit maps
$\rho_w$, for $w \in W$, are locally constant if and only if the
$\gfr$-action on $W$ is trivial or equivalently if and only if the
closed 2-sided ideal $I(\gfr )$ in $D(G,K)$ generated by $\gfr$
annihilates $W$.

What can we say about the quotient algebra $D^{\infty}(G,K) :=
D(G,K)/I(\gfr )$ ?

$D(G,K)$ is the strong dual of the space $C^{an}(G,K)$ of
$K$-valued locally analytic functions on $G$. Since the Dirac
distributions generate a dense subspace in $D(G,K)$ ([ST] 3.1) the
ideal $I(\gfr)$ is the orthogonal of the closed subspace in
$C^{an}(G,K)$ which is the simultaneous kernel of all linear forms
$\delta_g *\xfr *\delta_h$ with $\xfr \in \gfr$ and $g,h \in G$.
This is precisely the subspace of those functions in $C^{an}(G,K)$
which are annihilated by the action of $\gfr$. And this in turn,
by Taylor's formula, is the subspace $C^{\infty}(G,K)$ of all
$K$-valued locally constant functions on $G$ with the subspace
topology. On the other hand as a direct product of spaces of
compact type the space $C^{an}(G,K)$ is reflexive. In this
situation the strong dual of a closed subspace is the quotient of
the strong dual by the orthogonal subspace ([B-TVS] IV.16 Cor.).
In other words we have
$$
D^{\infty}(G,K) = C^{\infty}(G,K)'_b\ .
$$
Moreover, if $H \subseteq G$ is a fixed compact open subgroup,
then
$$
C^{\infty}(G,K) = \prod_{g\in G/H} C^{\infty}(gH,K)
$$
is the direct product of the spaces $C^{\infty}(gH,K)$ each of
which is a locally convex inductive limit of finite dimensional
spaces and hence carries the finest locally convex topology
(compare [ST] 1.2.i). In particular $C^{\infty}(H,K)$ is the
inductive limit
$$
C^{\infty}(H,K) =
\mathop{\lim}\limits_{\mathop{\longrightarrow}\limits_{N}}\; K[H/N]
$$
of the algebraic group rings $K[H/N]$ with $N$ running through the
open normal subgroups of $H$.

All of this applies to $V^c$ for any admissible $G$-representation
$V$. In parti-cular, $V^c$ as well as its strong dual $(V^c )'_b$ are
$D^{\infty}(G,K)$-modules. Clearly $M := (V^c)'_b$ is an analytic
$D(G,K)$-module if and only if $M$ is finitely generated as a
$D^{\infty}(H,K)$-module for some fixed but arbitrary choice of a
compact open subgroup $H \subseteq G$. This condition can be expressed
purely in terms of multiplicities as follows. Let $\widehat{H}$ denote
the set of isomorphism classes of all irreducible smooth
$H$-representations. Recall that any $\pi \in
\widehat{H}$ is finite dimensional. We let
$$
\mu (\pi ) := \hbox{\rm multiplicity\ of}\ \pi\ \hbox{\rm in}\
C^{\infty}(H,K)
$$
so that we have
$$
C^{\infty}(H,K) \cong \mathop{\bigoplus}\limits_{\pi \in
\widehat{H}}
\mu (\pi )\cdot \pi
$$
and
$$
D^{\infty}(H,K) \cong \mathop{\prod}\limits_{\pi \in \widehat{H}}
(\pi^{\ast})^{\times \mu (\pi )}\leqno{(3)}
$$
where $\pi^{\ast}$ denotes the contragredient of $\pi$. Any smooth
$G$-representation $V$ is semisimple as an $H$-representation.
Moreover $V$ is admissible if and only if the multiplicities
$$
\mu (\pi ,V) := \hbox{\rm multiplicity\ of}\ \pi\ \hbox{\rm in}\ V
$$
for any $\pi \in \widehat{H}$ are finite. We then have
$$
V \cong \mathop{\bigoplus}\limits_{\pi \in
\widehat{H}}
\mu (\pi,V)\cdot \pi
$$
and
$$
(V^c )'_b \cong \mathop{\prod}\limits_{\pi \in \widehat{H}}
(\pi^{\ast})^{\times \mu (\pi,V)}\leqno{(4)}
$$
as $D^{\infty}(H,K)$-modules.

\medskip

{\bf Definition:}

{\it A smooth $G$-representation is called strongly admissible if
there is a natural number $m$ such that
$$
\mu (\pi ,V) \leq m\cdot\mu (\pi)
$$
for any $\pi \in \widehat{H}$.}

\medskip

That the above definition does not depend on the particular choice
of $H$ can be seen as follows. Let $H_0 \subseteq H$ be a pair of
compact open subgroups in $G$. For any $\pi \in \widehat{H}$ and
$\sigma \in \widehat{H_0}$ let $\mu (\pi :\sigma )$ denote the
multiplicity of $\sigma$ in $\pi | H_0$. One easily checks that
$$
[H:H_0]\cdot\mu (\sigma ) = \sum_{\pi \in \widehat{H}} \mu (\pi
:\sigma )\cdot \mu (\pi )\ \ \hbox{\rm and}\ \ \mu (\pi ) =
\sum_{\sigma \in \widehat{H_0}} \mu (\pi :\sigma )\cdot \mu (\sigma )\ .
$$
Assuming that $\mu (\pi ,V) \leq m\cdot\mu (\pi)$, resp. $\mu
(\sigma ,V) \leq n\cdot\mu (\sigma)$, we compute
$$
\matrix{\mu (\sigma ,V) & = & \mathop{\sum}\limits_{\pi \in \widehat{H}}
\mu (\pi :\sigma )\cdot\mu (\pi ,V)\hfill\cr\cr & \leq &
m\cdot\mathop{\sum}\limits_{\pi \in \widehat{H}}
\mu (\pi :\sigma )\cdot\mu (\pi )\hfill\cr\cr & = & m\cdot [H:H_0 ]\cdot\mu
(\sigma )\ ,\hfill}
$$
resp.
$$
\matrix{\mu (\pi ,V) & \leq & \mathop{\sum}\limits_{\sigma \in \widehat{H_0}}
\mu (\pi :\sigma )\cdot\mu (\sigma ,V)\hfill\cr\cr & \leq &
n\cdot\mathop{\sum}\limits_{\sigma \in \widehat{H_0}}
\mu (\pi :\sigma )\cdot\mu (\sigma )\hfill\cr\cr & = & n\cdot\mu (\pi )\ .\hfill}
$$

\medskip

{\bf Proposition 2.1:}

{\it The functor $V \longmapsto (V^c )'_b$ is an (anti)equivalence
of categories between the category of all strongly admissible
$G$-representations and the category of all analytic
$D(G,K)$-modules which are annihilated by $I(\gfr )$.}

Proof: Comparing $(3)$ and $(4)$ it is obvious that $(V^c )'_b$ is
finitely generated as a $D^{\infty}(H,K)$-module if and only if $V$ is
strongly admissible. Hence the functor in question is well defined and
fully faithful. Moreover, if $M$ is an analytic $D(G,K)$-module
annihilated by $I(\gfr )$ then we have a topological surjection
$D^{\infty}(H,K)^r \dlongrightarrow M$ for some $r
\geq 1$. The dual embedding $M' \hookrightarrow C^{\infty}(H,K)^r$
shows that $V := M'_b$ carries the finest locally convex topology
and therefore is a strongly admissible $G$-representation. By
reflexivity we have $M = (V^c )'_b$ so that $M$ lies in the image
of our functor.

\medskip

{\bf Proposition 2.2:}

{\it If $G$ is the group of $L$-rational points of a connected
reductive $L$-group $\bf G$ then any smooth $G$-representation of
finite length is strongly admissible.}

Proof: Let $C$ be a fixed algebraically closed field which
contains $K$. We first want to reduce the assertion to the case
where the coefficient field of the smooth representation is $C$.
Denoting by $(.)_C$ the base extension functor from $K$ to $C$ we
have
$$
V_C \cong \mathop{\bigoplus}\limits_{\pi \in
\widehat{H}} \mu (\pi,V)\cdot \pi_C\ .
$$
Let ${\rm Irr}_C(H)$ denote the set of isomorphism classes of all
irreducible smooth $H$-representations over $C$. For each $\sigma
\in  {\rm Irr}_C(H)$ there is a unique $\pi(\sigma) \in
\widehat{H}$ such that $\sigma$ occurs in $\pi(\sigma)_C$. The
theory of the Schur index tells us the following ([CR] (70:15)):

1) The Schur index $m_K(\sigma)$ of $\sigma \in {\rm Irr}_C(H)$
with respect to $K$ only depends on $\pi(\sigma)$; we therefore
put $m_K(\pi) := m_K(\sigma)$ if $\pi = \pi(\sigma)$.

2) For any $\pi \in \widehat{H}$ we have the decomposition
$$
\pi_C \cong m_K(\pi)\cdot
\mathop{\oplus}\limits_{\pi(\sigma)=\pi} \sigma\ .
$$

3) If $\pi = \pi(\sigma)$ then $\mu(\pi)\cdot m_K(\pi) = {\rm dim}_C\,
\sigma$.

By using 2) our above decomposition of $V_C$ becomes
$$
V_C \cong \mathop{\bigoplus}\limits_{\sigma \in {\rm Irr}_C(H)}
\mu (\pi(\sigma),V)\cdot m_K(\pi(\sigma))\cdot\sigma\ .
$$
If we therefore show that there is an $m \in \Ndss$ such that
$\mu(\sigma,V_C) = \mu(\pi(\sigma),V)\cdot m_K(\pi(\sigma)) \leq
m\cdot {\rm dim}_C\,\sigma$ for any $\sigma \in {\rm Irr}_C(H)$ then
it follows from 3) that $\mu(\pi,V) \leq m\cdot\mu(\pi)$ for any $\pi
\in
\widehat{H}$. According to [Vig] II.4.3.c with $V$ also $V_C$ is of
finite length. This reduces us to proving our assertion for smooth
$G$-representation over some algebraically closed field $C$ containing
the field of complex numbers $\Cdss$. We first look at the case when
$V$ is irreducible supercuspidal. By a character twist we may assume
that the central character of $V$ is of finite order. According to
[Vig] II.4.9 the representation $V$ then is the base extension to $C$
of an irreducible supercuspidal $G$-representation over $\Cdss$. For
the latter our assertion is a theorem of Harish-Chandra ([HC] Cor. of
Thm. 2), and it is obvious that base extension between two
algebraically closed fields respects our assertion. Since a general
irreducible $V$ is contained in a representation parabolically induced
from a supercuspidal representation it remains to show that parabolic
induction respects strong admissibility. Let $P = P_LP_u$ be a
parabolic subgroup of $G$ with unipotent radical $P_u$ and Levi factor
$P_L$ and let $W$ be a strongly admissible smooth representation of
$P_L$. We have to check that $V := {\rm Ind}^G_{P}(W)$ is again
strongly admissible. Since $V$ is known to be admissible ([Vig]
II.2.1) we can do this by proving that the full linear dual $V^{\ast}$
of $V$ is finitely generated as an $D^{\infty}(H,K)$-module. Moreover,
being completely free in the choice of the compact open subgroup $H$
of $G$ we may choose it in such a way that the Iwasawa decomposition
$G = HP$ holds. Put $H_P := H \cap P$ and let $H_L$ denote the image
of $H_P$ in $P_L$. As an $H$-representation we then have
$$
{\rm Ind}^G_P(W) = {\rm Ind}^H_{H_P}(W|H_L)\ .
$$
By assumption $(W|H_L)^{\ast}$ is a finitely generated
$D^{\infty}(H_L,K)$-module. All we have to see therefore is that
$$
{\rm Ind}^H_{H_P}(W|H_L)^{\ast} = D^{\infty}(H,K)
\mathop{\otimes}\limits_{D^{\infty}(H_P,K)} (W|H_L)^{\ast}
$$
holds true. By semisimplicity this is an easy consequence of the
analogous identity with $W|H_L$ replaced by $C^{\infty}(H_P,K)$.

\medskip

As a consequence of these results we obtain that the functor $V
\longmapsto (V^c )'_b$ induces a bijective correspondence between
irreducible smooth $G$-representations and quasi-simple analytic
$D(G,K)$-modules which are annihilated by $I(\gfr )$. It should be
pointed out that $(V^c )'_b$ as a vector space is the full linear
dual of $V$. The smooth linear forms form a in general proper
$D^{\infty}(G,K)$-submodule of $(V^c )'_b$ so that the latter
cannot be simple.

\bigskip

{\bf 3. $U(\gfr)$-finite modules}

\smallskip

In this section we let $G$ be the group of $L$-rational points of
a connected reductive split $L$-group $\bf G$. We want to
understand more generally those analytic $D(G,K)$-modules $M$ on
which $U(\gfr)$ acts through a finite dimensional quotient. They
will be called $U(\gfr)$-finite.

Let $E$ be the underlying $L$-vector space of an irreducible
$L$-rational algebraic representation of $\bf G$. For any
$U(\gfr)$-finite analytic $D(G,K)$-module $M$ we set
$$
M^E := \Hom_{U(\gfr)}(E,M)\ .
$$
$\Hom_L (E,M)$ and hence $M^E$ as a closed vector subspace both
inherit a natural Fr\'echet topology from $M$. The group $G$ acts
on $M^E$ via the continuous $K$-linear endomorphisms
$$
f \longmapsto \,^g f(x) := g(f(g^{-1}x))\ \ \hbox{\rm for}\ g \in
G\ \hbox{and}\ f \in M^E .
$$
Moreover,
$$
\matrix{E\ \times\ M^E & \longrightarrow & M\cr (x,f) &
\longmapsto & f(x)}
$$
is a continuous $G$-equivariant bilinear map.

Let $V := M'_b$ denote the strong dual of $M$ as a locally
analytic $G$-representation. In order to determine the topology on
$V$ we need the following result.

\medskip

{\bf Proposition 3.1:}

{\it Let $J\subseteq U(\gfr)$ be a 2-sided ideal of finite
codimension and let $H\subseteq G$ be a compact open subgroup;
then the subspace topology on the subspace $C^{an}(H,K)^{J=0}$ of
all vectors in $C^{an}(H,K)$ annihilated by $J$ is the finest
locally convex topology.}

Proof: Fix an ordered vector space basis for $\gfr$, and an
exponential map for $G$.  This data, together with a choice of disk of
sufficiently small radius $s$ around the origin in $L^{{\rm
dim}\gfr}$, determines a ``canonical chart of the second kind'' on
$H$. Let $H_{r}$  be the family of standard compact open subgroups of
$H$ obtained from this canonical chart (see [Fea] 4.3.3). The Banach
space of analytic functions on $H_{r}$ is the standard Banach space
${\cal F}_{0,r}(K)$ of convergent series with coefficients in $K$ on
the disk of radius $r$ for $0<r\le s$. Let
$$
{\cal F}_{r}:=\prod_{h\in H_{r}\backslash H} {\cal F}_{0,r}(K)\
.
$$
Following the proof of [Fea] 3.3.4 we see that this Banach space
is an analytic $H_{s}$-representation and
$$
\lim\limits_{\to}{\cal F}_{r}\;\mathop{\longrightarrow}\limits^
{\sim}\;C^{an}(H,K)\ .
$$
By [Fea] 4.7.3, there is a non-degenerate pairing between
$U(\gfr)$ and the factor ${\cal F}_{0,r}$ of the product defining
${\cal F}_{r}$ corresponding to the trivial coset $H_{r}$.  This
pairing is given by evaluation at the identity element:
$$
\matrix{
U(\gfr)\times {\cal F}_{0,r} & \to & K\cr
 (\zfr,f) & \mapsto & (\zfr f)(1)\ .
}
$$
The ideal $J$ is of finite codimension in $U(\gfr)$, and given the
non-degeneracy of the pairing it follows that the space ${\cal
F}_{0,r}^{J=0}$ is finite dimensional.  Furthermore, because the
$U(\gfr)$-action from the left commutes with the right translation
action of $H$ it follows immediately that ${\cal F}_{r}^{J=0}$ is
finite dimensional.  Then $C^{an}(H,K)^{J=0}$, being the direct
limit of these finite dimensional spaces ([ST] 1.2.i), has the
finest locally convex topology.

\medskip

Since $M$ is analytic we have a surjection $D(H,K)^m
\dlongrightarrow M$ of $D(H,K)$-modules for some $m \in \Ndss$ and
some (or any) compact open subgroup $H \subseteq G$. After
dualizing we obtain an injection $V \hookrightarrow C^{an}(H,K)^m$
which certainly is $U(\gfr)$-linear ([ST] 3.2). Moreover, by
assumption there is a 2-sided ideal $J \subseteq U(\gfr)$ of
finite codimension which annihilates $M$ and hence $V$. Hence we
actually have an injection $V \hookrightarrow
(C^{an}(H,K)^{J=0})^m$. Applying Prop. 3.1 we now see that the
topology on $V$ necessarily is the finest locally convex one.

For general reasons $E \otimes_L V$ with $G$ acting diagonally
also is a locally analytic $G$-representation on a $K$-vector
space of compact type ([Fea] 2.4.3 and [ST] 1.2.ii). For our
particular $V$ the topology on $E \otimes_L V$ is, according to
the above discussion, the finest locally convex one. We let $E
\otimes_{U(\gfr)} V$ denote the $G$-equivariant quotient of $E \otimes_L V$
by the (automatically closed) $K$-vector subspace generated by all
vectors of the form $\xfr x \otimes v + x \otimes \xfr v$ for
$\xfr \in \gfr$, $x \in E$, and $v \in V$. By [ST] 1.2.i this
quotient is a locally analytic $G$-representation on a $K$-vector
space of compact type whose topology is the finest locally convex
one and whose strong dual evidently is $M^E$. In particular, both
$E \otimes_{U(\gfr)} V$ and $M^E$ are separately continuous
$D(G,K)$-modules.

By continuity and [ST] 3.1 the above bilinear map $E \times M^E
\longrightarrow M$ induces a continuous $D(G,K)$-module
homomorphism
$$
E \otimes_L M^E \longrightarrow M\ .
$$
By construction the $\gfr$-action on $E \otimes_{U(\gfr)} V$
derived from the $G$-action is trivial. Hence $I(\gfr )$
annihilates $E \otimes_{U(\gfr)} V$ and by duality also $M^E$.
Provided that $M^E$ is finitely generated as a $D(H,K)$-module for
some compact open subgroup $H \subseteq G$ it follows from Prop.
2.1 that $M^E$ is the dual of the strongly admissible
$G$-representation $E \otimes_{U(\gfr)} V$.

Let $\widehat{\bf G}$ denote the set of isomorphism classes of all
irreducible $L$-rational algebraic representations of $\bf G$. We
have the continuous $D(G,K)$-module module homomorphism
$$
\bigoplus_{E \in \widehat{\bf G}} E \otimes_L M^E \longrightarrow
M\ .
$$
The direct sum on the left hand side in fact is finite since the
number of $E \in \widehat{\bf G}$ which are annihilated by a given
2-sided ideal of finite codimension in $U(\gfr)$ is finite.

\medskip

{\bf Proposition 3.2:}

{\it Assume that $\bf G$ is split semisimple and simply connected;
for any $U(\gfr)$-finite analytic $D(G,K)$-module $M$ the natural
map
$$
\bigoplus_{E \in \widehat{\bf G}} E \otimes_L M^E
\mathop{\longrightarrow}\limits^{\cong} M
$$
is an isomorphism of $D(G,K)$-modules, each $M^E$ is the linear
dual of a strongly admissible $G$-representation over $K$, and
$M^E = 0$ for all but finitely many $E \in \widehat{\bf G}$.}

Proof: We have already noted that the map in question is a
homomorphism of $D(G,K)$-modules and that the direct sum on the left
hand side is finite. To establish the bijectivity we set $\gfr_K :=
\gfr \otimes_L K$ and we let $\widehat{\gfr}$, resp.
$\widehat{\gfr_K}$, denote the set of isomorphism classes of all
finite dimensional simple $\gfr$-modules, resp. $\gfr_K$-modules. By
assumption $M$ is an $U(\gfr_K)/J$-module for some 2-sided ideal $J
\subseteq U(\gfr_K)$ of finite codimension. Since $\gfr$
is semisimple the algebra $U(\gfr_K)/J$ is semisimple ([Dix]
1.6.4). Hence $M$ is a semisimple $\gfr_K$-module and we have its
isotypic decomposition
$$
M = \bigoplus_{E \in \widehat{\gfr_K}} M_E
$$
(compare [Dix] 1.2.8). Moreover, since ${\rm End}_{U(\gfr_K)}(E) =
K$ ([Dix] 2.6.5 and 7.2.2(i)) we have the natural isomorphism
$$
E \otimes_K \Hom_{U(\gfr_K)}(E,M)
\mathop{\longrightarrow}\limits^{\cong} M_E
$$
for any $E \subseteq \widehat{\gfr_K}$. Since the functor $E
\longmapsto E \otimes_L K$ induces a bijection $\widehat{\gfr}
\mathop{\longrightarrow}\limits^\sim \widehat{\gfr_K}$ (both sides
are classified by the dominant weights) the above isotypic
decomposition can be rewritten as a bijection
$$
\bigoplus_{E \in \widehat{\gfr}} E \otimes_L \Hom_{U(\gfr)}(E,M)
\mathop{\longrightarrow}\limits^{\simeq} M\ .
$$
But since $\bf G$ is assumed to be simply connected we have, by
derivation, the natural bijection $\widehat{\bf G}
\mathop{\longrightarrow}\limits^\sim \widehat{\gfr}$ so that the
last bijection coincides with the isomorphism in the assertion.

With $M$ also its direct summand $E \otimes_L M^E$ is finitely
generated as a $D(H,K)$-module for any compact open subgroup $H
\subseteq G$. It follows that $M^E$ is a finitely generated
$D(H,K)$-module as well: Take finitely many tensors which generate
$E \otimes_L M^E$; their $M^E$-components form a generating set
for $M^E$. We have explained above that then $M^E$ is the linear
dual of a strongly admissible $G$-representation.

\medskip

{\bf Example:} The assumptions on the group $\bf G$ in the above
Proposition cannot be weakened as the following example shows. Let
$L = K := \Qdss_2$, ${\bf G} := {\bf PGL_3}$, and ${\bf G_o} :=
{\bf SL_3}$. Then $G_{\rm o} := {\bf SL_3}(\Qdss_2 )$ is an open
normal subgroup of index three in $G = {\bf PGL_3}(\Qdss_2 )$. Let
$E_{\rm o}$ denote the three dimensional standard representation
of $G_{\rm o}$ and let $M := {\rm Ind}_{G_{\rm o}}^G (E_{\rm o})$
be the induced $G$-representation (in the sense of abstract
groups). It is clear that $M$ is an $U(\gfr)$-finite analytic
$D(G,K)$-module. One checks that as a $G_{\rm o}$-representation
$M$ is isomorphic to $E_{\rm o} \oplus E_{\rm o} \oplus E_{\rm
o}$. Since $\widehat{\bf G}$ is a subset of $\widehat{\bf G_o} =
\widehat{\gfr}$ to which $E_{\rm o}$ does not belong we see that
$\Hom_{U(\gfr)}(E,M) = 0$ for any $E \in \widehat{\bf G}$.

\medskip

For the sake of completeness we remark that vice versa any finite
direct sum $E_1 \otimes_L \Hom_K (V_1 ,K) \oplus \ldots
\oplus E_r \otimes_L \Hom_K (V_r ,K)$ with $E_i \in \widehat{\bf
G}$ and strongly admissible smooth $G$-representations $V_i$ over
$K$ is an $U(\gfr)$-finite analytic $D(G,K)$-module. Apart from
the finite generation which is contained in the subsequent lemma
this is clear.

\medskip

{\bf Lemma 3.3:}

{\it Let $H \subseteq G$ be a compact open subgroup; for any
finitely generated $D^{\infty}(H,K)$-module $N$ and any $E \in
\widehat{\bf G}$ the $D(H,K)$-module $E \otimes_L N$ is finitely
generated.}

Proof: We begin with a general observation. Let $\Oscr (\bf G)$ denote
the space of $L$-rational functions on $G$. Then the map
$$
\matrix{
\Oscr ({\bf G})\,\mathop{\otimes}\limits_L\,C^{\infty}(H,K) &
\longrightarrow & C^{an}(H,K)\cr (\psi,f) & \longmapsto & \psi|H\cdot
f\hfill }
$$
is injective. This can be seen as follows. Let $\sum_{j=1}^m \psi_j
\otimes f_j$ be an element in the left hand side such that $\sum_j
\psi_j|H\cdot f_j = 0$. We may assume that $\psi_1,\ldots,\psi_m$ are
linearly independent. Choose a disjoint covering $H =
\dot{\bigcup}_{i=1}^n\,U_i$ by nonempty open subsets $U_i \subseteq H$
such that the restrictions $f_j|U_i$, for any $1 \leq j \leq m$ and $1
\leq i \leq n$, are constant. Since each $U_i$ is Zariski dense in
$\bf G$ (this can be deduced, e.g., from [DG] II.5.4.3 and
II.6.2.1) it follows that $\psi_1|U_i,\ldots,\psi_m|U_i$ viewed in
$C^{an}(U_i,K)$ still are linearly independent. Hence $f_j|U_i =
0$ for any $i$ and $j$ and therefore $f_j = 0$ for any $j$.

Coming back to our assertion it suffices, of course, to consider the
case $N = D^{\infty}(H,K)$. On the other hand, if $J \subseteq
U(\gfr)$ denotes the annihilator ideal of $E^{\ast}$ then we find some
$G$-equivariant embedding $E^{\ast}\hookrightarrow \Oscr({\bf
G})^{J=0}$. Combining this with the above map leads, using the Leibniz
rule, to an $H$-equivariant embedding $E^{\ast} \otimes_L
C^{\infty}(H,K) \hookrightarrow C^{an}(H,K)^{J=0} \subseteq
C^{an}(H,K)$. As a consequence of Prop. 3.1 the topology induced by
$C^{an}(H,K)$ on the left hand side is the finest locally convex
topology. By dualizing we therefore obtain a surjection $D(H,K)
\dlongrightarrow E \otimes_L D^{\infty}(H,K)$ of $D(H,K)$-modules.

%It suffices, of course, to consider the case $N = D^{\infty}(H,K)$.
%Since $D^{\infty}(H,K) \cong
%\mathop{\prod}\limits_{\pi \in \widehat{H}} \pi^{\times
%\mu (\pi^{\ast})}$ we are further reduced (??) to proving the following
%more precise assertion: Fix an $L$-basis $x_1 ,\ldots,x_m$ of $E$;
%for any $\pi \in \widehat{H}$ and any nonzero vector $v \in \pi$
%the vectors $x_1 \otimes v,\ldots,x_m \otimes v$ generate the
%$D(H,K)$-module $E \otimes_L \pi$. In order to see this let $M$
%denote the $D(H,K)$-submodule of $E \otimes_L \pi$ generated by
%the vectors $x_1 \otimes v,\ldots,x_m \otimes v$. Clearly $M$
%contains $E \otimes v$. Consider the vector subspace $\pi_{\rm o}
%:= \{v' \in \pi : E \otimes v' \subseteq M\}$ of $\pi$. We claim
%that $\pi_{\rm o}$ is a $D^{\infty}(H,K)$-submodule: For $h \in H$
%and $v' \in \pi_{\rm o}$ we have
%$$
%E \otimes hv' = hE \otimes hv' = h(E \otimes v') \subseteq hM = M\
%.
%$$
%Since $v \in \pi_{\rm o}$ it follows that $\pi_{\rm o} = \pi$ and
%consequently that $M = E \otimes_L \pi$.

\medskip

We finally study the question when an $U(\gfr)$-finite analytic
$D(G,K)$-module is quasi-simple.

\medskip

{\bf Proposition 3.4:}

{\it If $E \in \widehat{\bf G}$ and $V$ is an irreducible smooth
$G$-representation over $K$ then $E \mathop{\otimes}\limits_L V$ with
the diagonal $G$-action is a simple module over the group ring
$K[G]$.}

Proof: We show that each nonzero element $x \in E
\otimes_L V$ generates $E
\otimes_L V$ as a $K[G]$-module. But first we recall a
few facts from rational representation theory (compare [Jan] II
\S\S1 and 2). Fix a Borel subgroup $P
\subseteq G$ and a maximal split torus $T \subseteq P$, and let $N$
denote the unipotent radical of $P$.

1. The subspace $E^N$ of $N$-invariants in $E$ is one dimensional and
coincides with the weight space $E_{\lambda}$ where $\lambda$ is the
highest weight of $E$ (w.r.t. $T$ and $B$).

2. If $e \in E_{\mu}$ has weight $\mu$ then $Ne \subseteq e +
\sum_{\mu < \nu} E_{\nu}$.

The fact 1. holds true similarly on the level of Lie algebras. This
shows that whenever $U_{\rm o} \subseteq N$ is an open subgroup then

1'. $E^{U_{\rm o}} = E^N = E_{\lambda}$ is one dimensional.

Since $E$ is also an irreducible module for the induced action of the
Lie algebra of $G$ it follows that whenever $U \subseteq G$ is an open
subgroup we have

3. $L[U]\cdot e = E$ for any nonzero $e \in E$.

Consider now a fixed nonzero element
$$
x = e_1 \otimes v_1 + \ldots + e_r \otimes v_r
$$
with $0 \neq e_i \in E$ and $0 \neq v_i \in V$. We may assume that
each $e_i$ is a weight vector. In order to show that $K[G]\cdot x = E
\otimes_L V$ we may replace $x$ when convenient by any other nonzero
element in $K[G]\cdot x$. In a first step we will show that for this
reason we may assume in fact that $r = 1$.

By the smoothness of $V$ we find an open subgroup $U \subseteq G$
which fixes each of the vectors $v_1 ,\ldots,v_r$. Put $U_{\rm o} := U
\cap N$. If $U_{\rm o}$ fixes $x$ we are immediately reduced to the
case $r = 1$ since, by 1'., we then have $x \in (E \otimes_L
V^U)^{U_{\rm o}} = E^{U_{\rm o}} \otimes_L V^U = E_{\lambda} \otimes_L
V^U$. Otherwise there is a $g \in U_{\rm o}$ such that $gx - x \neq 0$
and we replace $x$ by
$$
gx - x = (ge_1 - e_1 ) \otimes v_1 + \ldots + (ge_r - e_r ) \otimes
v_r\ .
$$
The point to note is that, by 2., each $ge_i - e_i$ lies in a sum of
weight spaces where the occuring weights are strictly bigger than the
weight of $e_i$. This means that one way or another after finitely
many steps we have replaced $x$ by a nonzero element in $E_{\lambda}
\otimes_L V^U$ for which $r$ can be assumed to be one.

Let therefore, for the second step of the proof, $x \in E
\otimes_L V$ be an element of the form $x = e \otimes v$ with $0
\neq e \in E$ and $0 \neq v \in V$. Denote by $U \subseteq G$ the
stabilizer of $v$. Using 3. and the irreducibility of $V$ we
obtain
$$
K[G]\cdot x = K[G]\cdot ((L[U]\cdot e)\otimes v) = K[G]\cdot (E
\otimes v) = E \otimes K[G]\cdot v = E \otimes V\ .
$$

\eject
%\medskip

{\bf Corollary 3.5:}

{\it Assume that $\bf G$ is split semisimple and simply connected
and let $M$ be any $U(\gfr)$-finite analytic $D(G,K)$-module; then
$M$ is quasi-simple if and only if it is of the form $M \cong E
\otimes_L \Hom_K (V,K)$ for some $E \in \widehat{\bf G}$ and some
irreducible smooth $G$-representation $V$ over $K$.}

Proof: If $E \in \widehat{\bf G}$ and $V$ is irreducible smooth
then $E^{\ast} \otimes_L V$ is a simple $D(G,K)$-module by Prop.
3.4. Hence $E \otimes_L \Hom_K (V,K) = (E^{\ast} \otimes_L V)'$ is
quasi-simple.

If on the other hand $M$ is quasi-simple then there is, by Prop.
3.2, an $E \in \widehat{\bf G}$ and a strongly admissible
$G$-representation $V$ such that $M = E \otimes_L \Hom_K (V,K)$.
With $M$ also $\Hom_K (V,K)$ is quasi-simple. Hence $V$ is
irreducible.

\medskip

The results of this section have more or less obvious counterparts
for $G$ being a compact open subgroup in ${\bf G}(L)$. We leave
precise formulations to the reader.

\bigskip

{\bf 4. An example}

\smallskip

In this last section we will analyze the reducible members of the
locally analytic principal series of the group ${\bf SL_2}(\dQ)$ and
we will show that they contain tensor product representations of the
kind considered in the last section.

Throughout this section let $G := {\bf SL_2}(\dQ)$. Furthermore, let
$P$ denote the Borel subgroup of lower triangular matrices in $G$ and
$T$ the subgroup of diagonal matrices. We actually will view $T$ as a
quotient of $P$. Assuming that $K$ is contained in the completion of
an algebraic closure of $\dQ$ we fix a $K$-valued locally analytic
character
$$
\chi:T\to K^{\times}\ .
$$
The corresponding principal series representation is
$$
\Ind_{P}^{G}(\chi) := \{f\in C^{an}(G,K): f(gp)=
\chi(p^{-1})f(g) \ \hbox{for any}\ g\in G, p\in P\}
$$
with $G$ acting by left translation. This is a locally analytic
$G$-representation on a vector space of compact type and its strong
dual
$$
M_{\chi}:=\Ind_{P}^{G}(\chi)'_{b}
$$
is a $D(G,K)$-module which is finitely generated, e.g., as a
$D(B,K)$-module where $B$ is the Iwahori subgroup of $G$ ([ST] \S\S5
and 6). By Prop. 1.1 the $D(G,K)$-module $M_{\chi}$ therefore is
analytic.

The basic numerical invariant of the character $\chi$ which governs
the irreducibility properties of $\Ind_P^G (\chi)$ is the number
$c(\chi)\in K$ defined by the expansion
$$
\chi(\pmatrix{t^{-1} & 0\cr 0 & t})=\exp(c(\chi)\log(t))
$$
for $t$ sufficiently close to $1$ in $\dZ$. It is shown in [ST] 6.1
that $M_{\chi}$ is a simple $D(G,K)$-module if $c(\chi)
\not\in -\Ndss_0$. We therefore assume for the rest of this section
that $m := -c(\chi) \in \Ndss_0$. According to [ST] 6.2 we then have a
nonzero homomorphism of $D(G,K)$-modules $M_{\chi'} \longrightarrow
M_{\chi}$ where $\chi' := \epsilon^{m+1}\chi$ and
$\epsilon(\pmatrix{t^{-1} & 0\cr 0 & t}):= t^2$ is the positive root
of $G$ with respect to $P$. Since $c(\chi') = m+2$ the module
$M_{\chi'}$ is simple and the above map consequently is injective. It
therefore remains to discuss the quotient module
$$
M_{\chi}^{loc} := M_{\chi}/M_{\chi'}
$$
which, of course, is finitely generated. On the other hand, the above
map is exhibited in the proof of [ST] 6.2 as the dual $I'$ of a
$G$-equivariant continuous linear map
$$
I : \Ind_{P}^{G}(\chi) \longrightarrow \Ind_{P}^{G}(\chi')
$$
whose actual construction we will recall further below. By the
argument in [ST] 3.5 the kernel of $I$ again is a locally analytic
$G$-representation on a vector space of compact type. We will see that
$I$ is a quotient map or equivalently that the image $I'(M_{\chi'})$
is closed in $M_{\chi}$. The module $M_{\chi}^{loc}$ therefore is the
continuous dual of the kernel of $I$ and in particular is analytic.

Write $\chi = \chi_{alg}\cdot\chi_{lc}$ where
$\chi_{alg}(\pmatrix{t^{-1} & 0\cr 0 & t}):= t^{-m}$ is a
$\dQ$-rational character and $\chi_{lc}$ is a $K$-valued locally
constant character of $T$. The character $\chi_{alg}$ is dominant for
the Borel subgroup $P^-$ opposite to $P$; hence the algebraic
induction $\ind_P^G(\chi_{alg})$ is the irreducible $\dQ$-rational
representation of highest weight $\chi_{alg}$ (w.r.t. $P^-$) of $G$
(compare [Jan] II.2 and II.8.23). On the other hand, since the
character $\chi_{lc}$ is locally constant we may form the smooth
induced $G$-representation
$$
\Ind_{P,\infty}^{G}(\chi_{lc}) := \{f\in C^{\infty}(G,K): f(gp)=
\chi_{lc}(p^{-1})f(g) \ \hbox{for any}\ g\in G, p\in P\}
$$
over $K$ with $G$ acting by left translation. It is known ([Vig]
II.5.13) to be a smooth $G$-representation of finite length which, by
Prop. 2.2, implies that it is strongly admissible. There is the
obvious $G$-equivariant linear map
$$
\matrix{
\tau : & \ind_P^G(\chi_{alg}) \mathop{\otimes}\limits_{\dQ}
\Ind_{P,\infty}^{G}(\chi_{lc}) & \longrightarrow & \Ind_P^G
(\chi)\cr & (\psi,f) & \longmapsto & \psi\cdot f\ . }
$$
We claim that
$$
0 \longrightarrow \ind_P^G(\chi_{alg}) \mathop{\otimes}\limits_{\dQ}
\Ind_{P,\infty}^{G}(\chi_{lc}) \mathop{\longrightarrow}\limits^{\tau}
\Ind_P^G(\chi) \mathop{\longrightarrow}\limits^{I} \Ind_P^G(\chi')
\longrightarrow 0\leqno{(\ast)}
$$
is an exact sequence of locally convex $K$-vector spaces (where the
left hand term carries the finest locally convex topology). This means
that it is exact as a sequence of vector spaces and that the maps
involved are strict. By dualizing and observing that the rational
representations of $G$ are selfdual this leads to the following
result.

\medskip

{\bf Proposition 4.1:}

{\it If $c(\chi) \in -\Ndss_0$ then the $D(G,K)$-module
$M_{\chi}^{loc}$ is analytic and $U(\gfr)$-finite and is isomorphic to
the tensor product of the $\dQ$-rational $G$-representation
$\ind_P^G(\chi_{alg})$ and the full $K$-linear dual of the smooth
representation $\Ind_{P,\infty}^{G}(\chi_{lc})$ of finite length. }

\medskip

We begin by recalling the construction of $I$ from [ST] 6.2. The group
$G$ acts on $ C^{an}(G,K)$ by left and by right translations. Both
actions derive into an action of the Lie algebra $\gfr =
\sfr\lfr_2(\dQ)$. Whereas the actions coming from left translation are
denoted, as usual, by $f \mapsto gf$ for $g\in G$ and $f \mapsto \xfr
f$ for $\xfr\in \gfr$ we write $f \mapsto \xfr_rf$ for the
$\gfr$-action derived from right translation. Then
$$
I(f) = (\ufr^-)^{1+m}_r f
$$
where $\ufr^{-}:=\pmatrix{0 & 1\cr 0&0\cr}\in \gfr$.

Corresponding to the decomposition $G = BP\;\dot{\cup}\;BwP$ where $B
\subseteq G$ is the Iwahori subgroup and $w:= \pmatrix{0 & -1\cr 1 &
0}$ the sequence $(\ast)$ is the direct sum of the sequences
$$
0 \longrightarrow \ind_P^G(\chi_{alg}) \mathop{\otimes}\limits_{\dQ}
\Ind_{P,\infty}^{BP}(\chi_{lc}) \mathop{\longrightarrow}\limits^{\tau}
\Ind_P^{BP}(\chi) \mathop{\longrightarrow}\limits^{I} \Ind_P^{BP}(\chi')
\longrightarrow 0
$$
and
$$
0 \longrightarrow \ind_P^G(\chi_{alg})
\mathop{\otimes}\limits_{\dQ}
\Ind_{P,\infty}^{BwP}(\chi_{lc}) \mathop{\longrightarrow}\limits^{\tau}
\Ind_P^{BwP}(\chi) \mathop{\longrightarrow}\limits^{I} \Ind_P^{BwP}(\chi')
\longrightarrow 0\ .
$$
The superscripts $BP$ and $BwP$ indicate the subspaces of those
functions in the induced representation which are supported in $BP$
and $BwP$, respectively. Both these sequences can be computed
explicitly as follows. Let $U$, resp. $U^-$, be the unipotent radical
of $P$, resp. $P^-$, and define $U_{\rm o} := U \cap B$ and $U^-_{\rm
o} := U^- \cap B$. Denoting by $u$, resp. $u^-$, the function on
$U_{\rm o}$, resp. $U_{\rm o}^-$, which sends a matrix to its left
lower, resp. right upper, entry we introduce the finite dimensional
$\dQ$-vector spaces $Pol^m(U_{\rm o})$ and $Pol^m(U_{\rm o}^-)$ of
polynomials of degree $\leq m$ in $u$ and $u^-$, respectively, with
coefficients in $\dQ$. By restricting, resp. translating by $w$ and
restricting, functions the above two sequences become isomorphic to
$$
0 \longrightarrow Pol^m(U_{\rm o}^-) \mathop{\otimes}\limits_{\dQ}
C^{\infty}(U_{\rm o}^-,K) \mathop{\longrightarrow}\limits^{\tau}
C^{an}(U_{\rm o}^-,K) \mathop{\longrightarrow}\limits^{({d\over
du^-})^{1+m}} C^{an}(U_{\rm o}^-,K) \longrightarrow 0
$$
and
$$
0 \longrightarrow Pol^m(U_{\rm o}) \mathop{\otimes}\limits_{\dQ}
C^{\infty}(U_{\rm o},K) \mathop{\longrightarrow}\limits^{\tau}
C^{an}(U_{\rm o},K) \mathop{\longrightarrow}\limits^{(-{d\over
du})^{1+m}} C^{an}(U_{\rm o},K) \longrightarrow 0\ .
$$
In these sequences the injectivity of the first map as well as the
exactness in the middle are obvious. By Prop. 3.1 the subspace
topology on the kernel of the second map is the finest locally convex
topology. The surjectivity and strictness of the second map can either
be checked directly or can be seen as a special case of the more
general statement in [Fea] 2.5.4. This finishes the proof of the
exactness of $(\ast)$.

\medskip

The smooth $G$-representation $\Ind_{P,\infty}^{G}(\chi_{lc})$ is
of length at most 2. More precisely one has ([GGP] p.173) that
$\Ind_{P,\infty}^{G}(\chi_{lc})$ is irreducible except in the
following cases:

A) $\chi_{lc} = 1$ is the trivial character. Then
$\Ind_{P,\infty}^{G}(1)$ contains the one dimensional trivial
representation on the subspace of constant functions. The
corresponding quotient is the socalled Steinberg representation which
is irreducible.

B) $\chi_{lc}$ is the character $\chi_{lc}(\pmatrix{t^{-1} & 0\cr 0 &
t})= |t|^2$ where $|\ |$ denotes the normalized absolute value of
$\dQ$. Then $\Ind_{P,\infty}^{G}(\chi_{lc})$ contains the Steinberg
representation and the corresponding quotient is the one dimensional
trivial representation.

C) $\chi_{lc}$ is of the form $\chi_{lc}(\pmatrix{t^{-1} & 0\cr 0 &
t})= |t|\cdot \delta(t)$ for some non-trivial quadratic character
$\delta : \dQ^{\times} \longrightarrow K^{\times}$. Then
$\Ind_{P,\infty}^{G}(\chi_{lc})$ either is irreducible (but not
absolutely irreducible) or is the direct sum of two infinite
dimensional non-equivalent irreducible $G$-representations.

If we combine this information with Prop. 4.1 and Cor. 3.5 we obtain a
complete list of the quasi-simple constituents of $M_{\chi}^{loc}$ up
to isomorphism. In particular, each of them is isomorphic to the
tensor product of $\ind_P^G(\chi_{alg})$ and the full $K$-linear dual
of one of the irreducible smooth representations in the above list. At
this point it should be mentioned that the length of a Jordan-H\"older
series for the kernel of $I$ on $\Ind_P^G(\chi)$ was already
determined in [Mor].

\eject
%\bigskip

{\bf References}

\parindent=23truept

\smallskip

\ref{[B-TVS]} Bourbaki, N.: Topological Vector Spaces.
Berlin-Heidelberg-New York: Springer 1987

\ref{[CR]} Curtis, C.W., Reiner, I.: Representation theory of
finite groups and associative algebras. New York-London: Wiley
1962

\ref{[DG]} Demazure, M., Gabriel, P.: Groupes Alg\'ebriques.
Amsterdam: North-Holland 1970

\ref{[Dix]} Dixmier, J.: Enveloping Algebras. AMS 1996

\ref{[Fea]} F\'eaux de Lacroix, C. T.: Einige Resultate \"uber die
topologischen Darstellungen $p$-adischer Liegruppen auf unendlich
dimensionalen Vektorr\"aumen \"uber einem $p$-adischen K\"orper.
Thesis, K\"oln 1997, Schriftenreihe Math. Inst. Univ. M\"unster,
3. Serie, Heft 23, pp. 1-111 (1999)

\ref{[GGP]} Gel'fand, I.M., Graev, M.I., Pyatetskii-Shapiro, I.I.:
Representation Theory and Automorphic Functions. San Diego:
Academic Press 1990

\ref{[HC]} Harish-Chandra, van Dijk: Harmonic Analysis on
Reductive $p$-adic Groups. Lect. Notes Math., vol. 162.
Berlin-Heidelberg-New York: Springer 1970

\ref{[Jan]} Jantzen, J.C.: Representations of Algebraic Groups.
Orlando: Academic Press 1987

\ref{[Mor]} Morita, Y.: Analytic Representations of $SL_2$ over a
$p$-Adic Number Field, III. In Automorphic Forms and Number Theory,
Adv. Studies Pure Math. 7, pp.185-222. Tokyo: Kinokuniya 1985

\ref{[ST]} Schneider, P., Teitelbaum, J.: Locally analytic
distributions and $p$-adic representation theory, with
applications to $GL_2$. Preprint 1999

\ref{[Vig]} Vigneras, M.-F.: Repr\'esentations $l$-modulaires d'un
groupe r\'eductifs $p$-adique avec $l \neq p$. Progress in Math.,
vol. 137. Boston-Basel-Stuttgart: Birkh\"auser 1996

\vfill\eject
%\bigskip

\parindent=0pt

Peter Schneider\hfill\break Mathematisches Institut\hfill\break
Westf\"alische Wilhelms-Universit\"at M\"unster\hfill\break
Einsteinstr. 62\hfill\break D-48149 M\"unster, Germany\hfill\break
pschnei@math.uni-muenster.de\hfill\break
http://www.uni-muenster.de/math/u/schneider\hfill\break

Jeremy Teitelbaum\hfill\break Department of Mathematics,
Statistics, and Computer Science (M/C 249)\hfill\break University
of Illinois at Chicago\hfill\break 851 S. Morgan St.\hfill\break
Chicago, IL 60607, USA\hfill\break jeremy@uic.edu\hfill\break
http://raphael.math.uic.edu/$\sim$jeremy\hfill\break

\end